\RequirePackage{ifpdf}
\ifpdf 
\documentclass[pdftex]{sigma-SAN}
\else
\documentclass{sigma-SAN}
\fi

\usepackage{amssymb}
 \usepackage{amsthm}
\usepackage{empheq}
\usepackage{fancyhdr}
\usepackage{color,soul}
\usepackage{datetime}
\date{\currenttime}
\input colordvi.sty
\usepackage{graphicx}

\newenvironment{prf}{\noindent\\[-2ex]{\bf Proof }}{\hfill $\Box$\\}
\newcommand{\cS}{{\cal S}}
\newcommand{\cU}{{\cal U}}

     \def\DD{\mathbb{D}}

     \def\Xh{{\mathbb{I}\hspace{-.2em} \rm  h\kern1pt}}
         \def\DH{{{\rm I}\hspace{-.15em}{\rm D}}}

\newcommand{\R}   {{\mbox{R\hskip-0.9em{}I \ }}}

\begin{document}

\newtheorem{thrm}{Theorem}
\newtheorem{prop}[thrm]{Proposition}
\newtheorem{cor}[thrm]{Corollary}

\theoremstyle{remark}
\newtheorem*{rem}{Remark}

\renewcommand{\PaperNumber}{***}

\FirstPageHeading

\ShortArticleName{Recursion Operators admitted by non-Abelian Burgers equations
}

\ArticleName{Recursion Operators admitted by non-Abelian \\ Burgers equations: Some Remarks}
\Author{Sandra CARILLO~(1 and 2),  Mauro LO  SCHIAVO~(1) and Cornelia SCHIEBOLD~(3 and 4)}

\AuthorNameForHeading{S.~Carillo, M.~Lo  Schiavo and C.~Schiebold}

\Address{(1)~Dipartimento ``Scienze di Base e Applicate
    per l'Ingegneria'',    
    
    \hskip-.5em   \textsc{Sapienza} - Universit\`a di Roma,  16, Via A. Scarpa, 00161 Rome, Italy
        } 
\EmailD{\href{mailto:email@address}{sandra.carillo\symbol{64}sbai.uniroma1.it}} 
\URLaddressD{\url{http://www.sbai.uniroma1.it/~sandra.carillo}} 
\EmailD{\href{mailto:email@address}{mauro.loschiavo\symbol{64}sbai.uniroma1.it}} 
\URLaddressD{\url{http://www.sbai.uniroma1.it/~mauro.loschiavo}} 

\Address{(2)~I.N.F.N. - Sez. Roma1,
Gr. IV - 
Mathematical Methods in NonLinear Physics,  Rome, Italy      

        } 

\Address{(3)~Department of  Science Education and Mathematics
                           
            \hskip-.5em  Mid Sweden University, S-851 70 Sundsvall, Sweden} 
\EmailD{\href{mailto:email@address}{Cornelia.Schiebold\symbol{64}miun.se}} 
\URLaddressD{\url{http://www.miun.se/personal/corneliaschiebold}} 

\Address{(4)~Instytut Matematyki, Uniwersytet Jana Kochanowskiego w Kielcach, Poland}
\ArticleDates{Received ???, in final form ????; Published online ????}

\bigskip

\Abstract{The recursion operators admitted by different operator Burgers equations, in the
framework of the study of nonlinear evolution equations, are here considered.
Specifically,  evolution equations wherein the unknown is an operator acting on a 
Banach space
are investigated.  Here, the {\it mirror} non-Abelian Burgers equation is
considered: it can be written as $r_t = r_{xx} + 2r_x r$. The structural properties of the obtained 
recursion operator are studied; thus, it is proved to be a strong symmetry for the {\it mirror} non-Abelian Burgers equation as well as to be the hereditary.  
{{These results are proved via direct computations as well as via computer assisted manipulations; ad hoc routines are needed to treat non-Abelian  quantities and relations among them.}}
The obtained recursion operator  generates the {\it mirror} 
non-Abelian Burgers hierarchy.   The latter, when the unknown operator $r$ is replaced by a real valued function reduces to the usual (commutative) Burgers hierarchy. Accordingly, also the  recursion operator reduces to the usual Burgers one.}

%
%
%


\section{Introduction }
\label{S:1} Non-Abelian Burgers equations are here studied. The idea is to construct
different non-commutative counterparts of the Burgers equation in a real valued unknown.
Indeed,  the  non-Abelian Burgers (or non-commutative Burgers) as  equation
usually considered takes the form of the corresponding nonlinear evolution equation,
namely $s_t = s_{xx} + 2ss_x$. Here, the {\it mirror} non-Abelian Burgers equation is
considered: it can be written as $r_t = r_{xx} + 2r_x r$. Both these non-Abelian Burgers
equation are studied by Kupershmidt in \cite{Ku} who constructed the whole hierarchies
they generate in the case of matrix equations. Notably, here the unknown in the equations
under investigation are supposed to be operators on a suitable Banach space. Hence, these
unknown cannot be represented via finite dimensional matrices. More precisely, $r(x,t)$ is a 
bounded linear endomorphism on some Banach space. In applications, choices for the 
underlying Banach space include sequence spaces  and $L_2(\R)$, see 
\cite{Carl Schiebold}, \cite{Carl Schiebold 1}, \cite{Schiebold2010}. On the other hand, 
the results on
non-commutative hierarchies of finite dimensional matrix equations are naturally
included, as a particular case,  in the present study. Non-Abelian generalization of Burgers 
equation where the unknown are  finite dimensional matrices are constructed in Bruschi, Levi and
Ragnisco \cite{LRB}.

The present study is concerned about structural properties of non-Abelian Burgers
equations and represents a continuation of the results in  \cite{SIMAI2008,
Carillo:Schiebold:JNMP2012}, where the recursion operator of the non-Abelian Burgers
equation is obtained via a Cole-Hopf \cite{Cole:1951, Hopf:1950} transformation
linking the non-commutative heat equation and the non-commutative Burgers equation.
Then, the obtained  operator is proved to satisfy all the  required algebraic
properties to be the hereditary recursion operator which generates the non Abelian Burgers 
hierarchy. Specifically, it is both a strong symmetry and a hereditary operator. 
following the same approach, the mirror hierarchy is generated. Notably, it coincides
with the mirror hierarchy proposed by Kupershmidt in \cite{Ku}, who constructs a
recursive definition of the hierarchies. Here, the hierarchy  is recovered on application
of the Cole-Hopf transformation viewed as a particular case of B\"acklund transformation
and, hence,  the results by Fuchssteiner \cite{Fuchssteiner1979} and  Fokas and Fuchssteiner \cite{FokasFuchssteiner:1981} referring to B\"acklund transformations and
recursion operators can be applied. In particular, the recursion operator of the mirror
Non Abelian Burgers equations is obtained, combining the non-commutative Cole-Hopf
transformation with the trivial recursion operator admitted by the non-commutative linear heat equation.

The hierarchy of non-commutative Burgers equations (therein termed right-handed) as well
as the corresponding recursion operator in  \cite{Carillo:Schiebold:JNMP2012}  were,
independently, obtained by G\"urses, Karasu and Turhan \cite{GKT} on application of a
method, in \cite{GKS}, based on the Lax pair formulation. 

It should be mentioned that the present investigation is
part of a wide reasearch program which takes its origins in the study of structural
properties of nonlinear evolution equations, where the unknown is a real valued function,
and their connection with B\"acklund transformations \cite{{Rogers:Carillo:1987b}, BS1, Fuchssteiner:Carillo:1989a}.
In particular, this work continues the study, currently under further development, on
non-Abelian nonlinear evolution equations in \cite{[77]} -  \cite{CMS-2015}, \cite{Fuchssteiner:Carillo:1989a}, \cite{Schiebold 1}- \cite{Schiebold2010}.

\noindent
The material is organized as follows.
The opening Section 2 concerns  the 
{\it mirror non-Abelian Burgers} equation, termed also mirror non-commutative Burgers
equation. This equation is linked, via a {\it mirror} Cole-Hopf transformation to the
noncommutative heat equation. The corresponding  hierarchy is generated via subsequent
applications of the admitted recursion operator, denoted as $\Phi(r)$ which is later shown to be hereditary. 
Then, all the equations belonging to the {\it mirror} Burgers hierarchy follow  on 
subsequent applications of the operator $\Phi(r)$. Notably, this  {\it mirror} Burgers 
hierarchy is the same obtained  by Kupershmidt \cite{Ku}.

In the subsequent Section 3, the obtained operator $\Phi(r)$,  is proved to represent 
a strong symmetry admitted by the mirror non-Abelian Burgers equation. 

Sections 4 is devoted 
the hereditarines of the recursion operator $\Phi(r)$. 
Notably, there are different ways to prove the hereditariness of the recursion operator 
$\Phi(r)$. Indeed, as already pointed out, its construction via the Cole-Hopf transformation 
which links the mirror  non-Abelian Burgers equation to the non-commutative linear heat 
equation indicates it inherits 
such a property. Furthermore, the result can be proved via a 
direct computation, following  the lines of the proof given in 
\cite{Carillo:Schiebold:JNMP2012} where the hereditariness  of the recursion operator 
admitted by the usual non-Abelian Burgers equation is shown.
In addition, the proof can be constructed via a computer assisted method: this is presented in
Section 5 where the difficulties which arise when a computer algebra language is used in 
dealing with non-commutative quantities is pointed out.

Finally, an Appendix  devoted to a brief summary on the connection between recursion
operators and B\"acklund transformations and, in particular, to summarize those results
needed throughout the other Sections, closes this work.  In addition, a brief overview on
results previouly obtained on  the noncommutative  Burgers equation, its recursion
operator and the link with the heat equation are also recalled to help the reader.

\section{The \textbf{\textit{mirror non-Abelian Burgers}} hierarchy}
\label{S:2} In this Section, the  non-Abelian Cole-Hopf transformation $r=u_xu^{-1}$  is
applied to the heat equation to obtain a {\it mirror}  non-Abelian Burgers equation, as,
according to \cite{Ku}, we term it.
That is, consider the B\"acklund transformation:
\begin{equation}\label{B2}
B(u,r)=0~, ~~~ \text{where}~~~ B(u,r)= ru-u_x~
\end{equation}
which links the heat equation $u_t= u_{xx}$ to the {\it mirror}  non-Abelian Burgers
equation
\begin{equation}\label{mbu}
r_t \ = \ r_{xx}+2r_xr,
\end{equation}
where,
following the method in \cite{Carillo:Schiebold:JNMP2012}, $\ r_{xx}+2r_xr=\Phi(r) r_x$, when
$\Phi(r)$ denotes
the recursion operator admitted by  the mirror  non-Abelian Burgers equation. 
\begin{prop}
The operator $\Phi(r)$ is given by
\begin{equation}\label{mrec}
\Phi(r)= (D-C_r)(D+R_r)(D-C_r)^{-1}, ~~~ \text{where}~~~C_r:=\big[r, \cdot\big],
\end{equation}
i.e. $C_r$ denotes the commutator with $r$, and $R_{r}$ is the right multiplication by
${r}$.
\end{prop}
\begin{prf}
Given the Cole-Hopf transformation $B$ in \eqref{B2}, its directional derivatives are:
   \begin{eqnarray*}
       B_u[V] &=& \frac{\partial}{\partial \epsilon}\Big|_{\epsilon=0}
                                \Big( r(u+\epsilon V)-(u+\epsilon V)_x \Big)
              = rV-V_x, \\
       B_r[W] &=& \frac{\partial}{\partial \epsilon}\Big|_{\epsilon=0}
                                \Big( (r+\epsilon W)u-u_x \Big)
              = Wu,
    \end{eqnarray*}
namely, for $V\in T_u\cU$, $W\in T_r\cS$, it follows $B_u = L_r-D $ and $B_r=R_u$,
hence the transformation operator $T=-B_r^{-1}B_u$. Then, when $L_{r}$ denotes the left 
multiplication by ${r}$, the following identities
    $$\begin{aligned}
 DR_u&=R_uD+R_uR_r = R_u(D+R_r) \\
L_uDL_{u^{-1}}&=L_u\left(L_{u^{-1}}D-L_{u^{-1}}L_{u_x}L_{u^{-1}}\right) = (D-L_r) \\
  (D-L_r)R_u&=(R_uD+R_uR_r)-L_rR_u = R_u(D-C_r)\\
L_u(D+R_r)&=DL_u-L_rL_u+L_uR_r = (D-C_r)L_u \\
R_{u^{-1}} D R_u&= D+R_{u^{-1}}R_{u_x}=(D+R_r) \\
\end{aligned} $$
allow to write the transformation operator $T$ in the form
\begin{equation}
T=(D-C_r) R_{u^{-1}}; 
\end{equation}
Then, the recursion operator $\Phi(r)$, given in \eqref{mrec}, is obtained via
\[
\Phi=T D  T ^{-1},
\]
where $D$ is the {\it trivial} recursion operator admitted by the linear heat equation.
\end{prf}
Hence, the mirror non-Abelian Burgers hierarchy is represented by
\begin{equation} \label{mbu hierarchy}
   r_{t_n} = \Phi(r)^{n-1}r_x, \quad n\geq 1,
\end{equation}
the lowest members of which read
\begin{equation}\label{mbuh}
\begin{array}{ccc}
   r_{t_1} &=& r_x, \hfill \\
   r_{t_2} &=& r_{xx} + 2r_xr, \hfill \\
   r_{t_3} &=& r_{xxx} + 3r_{xx} r+ 3r_x^2 + 3 r_x r^2.
\end{array}
\end{equation}
Note that all the members of this hierarchy are obtained from the corresponding ones in
the non-Abelian Burgers hierarchy when left multiplication is replaced with right
multiplication. Furthermore, also in this case, the whole hierarchy is linked via a
Cole-Hopf mirror transformation, which now is \eqref{B2}, instead of $B(u,s)=us-u_x$.
Transformation \eqref{B2} connects corresponding members in the heat hierarchy
\eqref{heat hierarchy} to corresponding ones in the non-Abelian Burgers mirror hierarchy
\eqref{mbu hierarchy}.
 
{{The next Sections are devoted to stated and   prove the main  Theorem on properties of the operator $\Phi(r)$.}}

\section{The non-Abelian mirror Burgers recursion operator}

This Section is devoted to the operator $\Phi(r)$ and, in particular,  the
following theorem is the main result.

\begin{thrm}\label{T1}
   The  operator   given in (\ref{mrec}), i.e.
$$
\Phi(r)= (D-C_r)(D+R_r)(D-C_r)^{-1}
$$
represents the hereditary  recursion operator of the non-Abelian mirror Burgers equation.
\end{thrm}
\noindent To prove the Theorem  \ref{T1} the following steps are needed
\begin{itemize}
\item prove that the  operator $\Phi(r)$ 
is a strong symmetry for the base member hierarchy, i.e. $r_t=H(r)$, where $H(r)=r_x$;
\item prove that the  operator $\Phi(r)$ 
 is hereditary.

\end{itemize}
Then, combination of the two steps completes the proof since it allows to conclude, as in
\cite{Carillo:Schiebold:JMP2011}, that the  operator $\Phi(r)$ is  hereditary.  Then, 
according to \cite{Carillo:Schiebold:JMP2011}, it is a strong symmetry for 
all the higher order nonlinear evolution equations of the non-Abelian mirror Burgers
hierarchy \eqref{mbu hierarchy}. 

\medskip\noindent
\begin{prf} (of Theorem \ref{T1})
{\bf Step 1} is represented by the following
\begin{prop} \label{SS}   The  operator  $\Phi(r)$ 
is a strong symmetry for $r_t=H(r)$, where $H(r)=r_x$.
\end{prop}
\begin{prf} (of Proposition \ref{SS})
The proposition is proved   when \footnote{see the Appendix and \cite{Fuchssteiner1979}.} the condition
\begin{equation} \label{strong symmetry RE1}
   \Phi'(r)[H(r)] = [H',\Phi(r)]
\end{equation} is shown to hold.

First of all, note that since $H(r)=r_x$, then  $H'(r)=D$; thus, on substitution of both of them   the
relation to prove becomes
 \begin{equation} \label{strong symmetry RE}
   \Phi'(r)[r_x] = [D,\Phi(r)].
\end{equation}
For computational convenience, the operator $\Phi(r)$ is re-written in the equivalent form
\begin{equation}\label{mrec2}
\Phi(r)=D+R_r+L_{r_x}(D-C_r)^{-1}~
\end{equation}
where, respectively, $C_r,\, R_r$ and $L_r$ denote the commutator, right and left multiplication by $r$,
that is
\begin{equation*}
C_r \sigma:= [r, \sigma], ~ R_r\sigma:= \sigma r, ~~  L_r\sigma:=r \sigma, ~~\forall \sigma.
\end{equation*}
Direct computation proves the thesis. The  Fr\'echet derivatives of the operator
$\Phi(r)$, in \eqref{mrec2}, is
\begin{equation} \label{mBurgers operator derivative}
   \Phi'(r) [V] = R_V + L_{V_x}(D-C_r)^{-1} + L_{r_x} (D-C_r)^{-1} C_V (D-C_r)^{-1}.
\end{equation}
The latter follows since $\forall V$,  $C_r'[V]=C_V$, $R_r'[V]=R_V$,
$L_{r_x}'[V]=L_{V_x}$, and product rule is applied so that  the Fr\'echet derivative of
$(D-C_r)^{-1}$ follows \footnote{Recall that $ \big(\Gamma^{-1}(r)\big)'[V] =
\Gamma(r)^{-1} \big( - \Gamma'(r)[V] \big) \Gamma(r)^{-1}$ holds for an operator-valued
function $\Gamma(r)$. }
\begin{equation}
   ((D-C_r)^{-1})'[V] =  (D-C_r)^{-1} C_V (D-C_r)^{-1}.
\end{equation}
To evaluate $\Phi'(r)[r_x]$, let $V=r_x$ in  \eqref{mBurgers operator derivative},
\begin{equation} \label{1}
   \Phi'(r) [r_x] = R_{r_x} + L_{r_{xx}}(D-C_r)^{-1} + L_{r_x} (D-C_r)^{-1} C_{r_x} (D-C_r)^{-1}.
\end{equation}
Now, since $[D-C_r, D] = [D,C_r] = C_{r_x}$ implies $[D,(D-C_r)^{-1}] = (D-C_r)^{-1}
C_{r_x} (D-C_r)^{-1}$,  the right hand side gives
\begin{eqnarray} \lefteqn{
   [D, \Phi(r)] \ =\ [D, R_r] + [D,L_{r_x}(D-C_r)^{-1}] } \nonumber\\ \nonumber
       &=& R_{r_x} + L_{r_{xx}}(D-C_r)^{-1} + L_{r_x}[D,(D-C_r)^{-1}] \\ \label{auxx}
       &=& R_{r_x} + L_{r_{xx}}(D-C_r)^{-1} + L_{r_x} (D-C_r)^{-1} C_{r_x} (D-C_r)^{-1},
\end{eqnarray}
Comparison of (\ref{1}) with (\ref{auxx}) shows (\ref{strong symmetry RE}) and completes
the proof.
\end{prf}

\medskip\noindent
The next {\bf Step 2} needed to prove Theorem \ref{T1} is represented by the proof that the operator
$\Phi(r)$ is hereditary: this result is established in the next Section.

{\bf Remark} A computer algebra program (using a symbolic language) was constructed to
provide a computer assisted proof of the recursivity of the operator $\Phi(r)$ and the
hereditariness of the same operator. Note that one of the main difficulties to overcome
writing computer routines that may prove results concerning non-Abelian properties is
that in the symbolic language, by default, all the variables are assumed to commute.
Hence, non commutativity requires non trivial {\it ad hoc} routines.

On the other hand, in devising the computer assisted proof there is no need of
introducing the notion of equivalence between operators, a relation useful to simplify
the computations done by hand. For instance, in \cite{Carillo:Schiebold:JNMP2012},
equivalence relations are introduced to avoid the explicit  computation of those terms
whose contribution satisfies the due symmetry requirement. The computer algebra routines
we prepared straightly produces all the terms and, then,  verify symmetry after the
exchange of the two arbitrary fields therein and the consequent sum. Some of the details
are given in  Section 5.
\end{prf}

\section{The hereditariness of the non-commutative mirror Burgers recursion operator}

This section is devoted   to  the hereditariness of  non-commutative
mirror Burgers recursion operator. The definition of hereditariness, introduced in
\cite{Fuchssteiner1979} 
 in the context of nonlinear evolution equations,
represents a key tool since, a {{ strong symmetry}} (recursion operator according
to \cite{Olver}) which is also hereditary  represents a  strong symmetry also for each one the 
nonlinear evolution equations of the hierarchy, in this case \eqref{mbu hierarchy},  it generates.
That is, the property is inherited from one equation to the next one in the hierarchy  and, 
hence, to the whole hierarchy. 
{{Hereditariness (see the  definition in the Appendix) is an algebraic property: it can be verified 
when a bilinear form is checked to be symmetric with respect to the exchange between each other 
of two arbitrary chosen fields it acts on. }}

{\begin{thrm}\label{T2} {\bf (Statement)}
   The non-commutative mirror Burgers recursion operator given in (\ref{mrec}) is hereditary.
\end{thrm}}
\noindent The thesis of this Theorem  can be proved in various different ways.
\begin{enumerate}
\item Indeed, Fokas and Fuchssteiner \cite{FokasFuchssteiner:1981} proved that hereditary operators are mapped to
hereditary operators via B\"acklund transformations. This result can be applied to the
 non-commutative mirror Burgers recursion operator since it is obtained via the
 Cole-Hopf transformation of the trivial recursion operator $D$ admitted by the heat
 equation.
\item   Following the method  in \cite{Carillo:Schiebold:JNMP2012}, a direct proof
can be  constructed computing all the terms in \eqref{hereditary}. Note that, the notion of equivalence can be introduced to simplify the required computations. 
\item  In addition, via an {\it ad hoc} computer algebra program which verifies that \eqref{hereditary} holds true.
\end{enumerate}

This third choice is examined in the next Section. Note that the idea to employ computer algebra routines to investigate properties of recursion operators in not new,
see  \cite{FuOWi} for early results,  
and \cite{Bill2010} (and references therein) for recent developments on the subject.
However, all of them are concerned about non linear evolution equations where the unknown is a real valued function and hence, the devised routines, in different symbolic languages, are in an Abelian framework while the present investigation concerns non-Abelian operator unknowns.
\section{Computer assisted results}

To ease-up the proof of some of the analytic properties of recursion operators, 
a computer
algebra program (using a symbolic language) was constructed that provides an automatic
assisted achievement of the necessary steps. At first, proof of the recursivity of the
operator $\Phi(r)$ has been produced. Then, to also prove hereditariness of the same
operator a second computer algebra program has been realized. Clearly, computer algebra
is convenient when long and tedious computations are necessary, however it must be
noticed that other technical problems arise. For instance, one of the main difficulties
to overcome has been that of writing routines that proved results concerning non-Abelian
computations. Indeed in the symbolic language, by default, all the variables are assumed
to commute, and all the operations such as multiplications, derivatives, and similar, are
commutative by default. Hence, non commutativity required non trivial {\it ad hoc}
routines.

Specifically, automatic  proofs procedure developed along the following subsequent steps.

The first step concerned realizing that operator $\Phi(r)$ may be easily rewritten if a
convenient {\it derivation} is introduced, namely, let us introduce the operator: 
\hbox{$\DH:=(D-C_r)$.}
Indeed,  $\DH$ has all necessary and characteristic properties of a derivative
(linearity, Leibnitz rule, etc.), and in the course of computation it may be (and has
been) used and interpreted as a normal derivative, provided that its real meaning is kept
into account. This is not only to say that: when $\DH f(x(\eta))$ needs to be computed
then the result is $f_x-\big[ r,f\big]$, but also that, when any other algebraic rule is
concerned, the new derivative $\DH=:\frac{\partial}{\partial\eta} $ may replace the
former $\frac{\partial}{\partial x} \ = D$ derivative, until the   variable $x\in\R$ is
replaced back at its place.

The second step is then that of writing the operator $\Phi(r)$ by use of this new
convenient derivative $\DH$. Its consequent compact form, from (\ref{mrec}), is 
easily found to be: 
$$ \Phi(r)=\DH\left(\DH+L_r\right)\DH^{-1} 
$$
and since it clearly is $\DH r\equiv Dr$, then this compact form for $\Phi(r)$
immediately shows that the equation's hierarchy is simply given by
$$ \Phi^n(r) Dr=\DH(\DH+L_r)^n\ r \ . $$

In particular, the compact form for the {\it mirror} Burgers equation has the easy aspect
$$ r_t=\Phi(r)Dr=\DH(\DH+L_r)r \ = \ \DH^2r+\DH r^2 \ . $$

Third step has been that of confirming the recursivity property of $\Phi(r)$ by automatic
computation with use of this new operator $\DH$. To achieve this, its Fr\'echet
derivative is needed, yet obviously keeping in mind that $\DH$ is still a function of the
equation variable $r$, and hence that the following hold
$$ \left\{\begin{aligned}
\big( \DH \big)^{\prime }[V] &=(D-C_r)^{\prime }[V] \ = \ -C_V \\
\big( \DH^{-1} \big)^{\prime }[V]&=-\DH^{-1}\big(\DH\big)^{\prime }\DH^{-1} \ = \
(D-C_r)^{-1}C_V(D-C_r)^{-1} \ .\end{aligned}\right. $$ Consequently, the Fr\'echet
derivative of $\Phi(r)=D+R_r+L_{r_x}(D-C_r)^{-1}$, given in \eqref{mBurgers operator derivative}, that  is
$$ \Phi^{\prime }[V]= R_V+L_{V_x}(D-C_r)^{-1}+L_{r_x}\big( (D-C_r)^{-1}\big)^{\prime }[V] ~,$$
turns out to acquire the computational more convenient form
\begin{align}
\Phi^{\prime }[V]&= -C_V+L_V+L_{V_x}\DH^{-1}+L_{r_x}\DH^{-1}C_V\DH^{-1} \\
 &= R_V+L_{\DH V}\DH^{-1}+L_{[r,V]}\DH^{-1}+L_{r_x}\DH^{-1}C_V\DH^{-1} \ , \label{mm1}
\end{align}
where it must be recalled that the field $\DH V$ is in fact
$\frac{\partial}{\partial\eta} \ V= \frac{\partial}{\partial x} \ V-C_rV$.

Next step is that of the technical (long and tedious) computations of the desired
properties. The first one, recursivity, is first performed using the base member of the
hierarchy, according to which the condition $\Phi^{\prime }[r_x]-\big[ D,\Phi \big]=0$ is
verified. To prove this fact, the operator $\DH$ may be used as the (unique) derivative
operator with respect to the {\it new}  variable $\eta$, however it has still been kept in
mind that this is possible only by replacing the {\it old}  derivative
$\frac{\partial}{\partial x} = D $ by the operator $(\DH+C_r)$, and by using the
Fr\'echet derivative of $\Phi(r)$ with its form \eqref{mm1}. This is actually what it has
been done to confirm the explicit direct proof that is also provided in the previous Section.
Furthermore, also to check the automatic procedure, the  next hierarchy
member has been obtained:
\begin{equation} \label{mm2}
\Phi^{\prime }[H(r)]=\Big[ H^{\prime }(r),\ \DH(\DH+L_r)\DH^{-1} \Big]
\end{equation}
where $H(r)$ is the ({{symmetric}}) Burgers equation: $H(r)=\DH^2 r+\DH r^2$, and $H^\prime
$ is its $r-$derivative: $ H^\prime (r)=D^2+2R_r D+2L_{r_x} $ expressed in the new
coordinates (and remember that $r_\eta\equiv r_x$):
$$ H^\prime (r) = \DH^2-R_{r_\eta}+3L_{r_\eta}+L_{r^2}-R_{r^2}+2L_r\DH \ . $$

It is useful to remark here that, although the variable $\eta$ coincides with the
variable $x$, all the same, due to non-commutative asset, their two derivations
$\DH$ and $D$ are different, and may coincide only in the commutative case. Only for
convenience, we write here the common value of \eqref{mm2}:

$$\left( \begin{array}{c}
r_{\eta \eta}\\
rr_\eta\\ R_{r_{\eta\eta\eta}}\DH^{-1}\\
r_\eta r \\
2 R_{r}R_{r_{\eta\eta}}\DH^{-1}\\
2 R_{r_\eta}R_{r_{\eta}}\DH^{-1}\\
 R_{r}R_rR_{r_{\eta}}\DH^{-1}\\
\end{array}\right) \quad +\quad
\left( \begin{array}{c}
R_{r_{\eta}}\DH^{-1}R_{r_{\eta\eta}}\DH^{-1}\\
-R_{r_\eta} R_{r}R_{r}\DH^{-1}\\
-R_{r_{\eta}}\DH^{-1}r_{\eta\eta}\DH^{-1}\\
R_{r_\eta}\DH^{-1}R_{r}R_{r_{\eta}}\DH^{-1}\\
R_{r_\eta}\DH^{-1}R_{r_\eta}R_{r}\DH^{-1}\\
-R_{r_\eta}\DH^{-1}rr_\eta\DH^{-1}\\
-R_{r_\eta}\DH^{-1}r_\eta r\DH^{-1}\\
\end{array}\right) $$

It may also be remarked that, although this being only a matter of chance, the same
result may be found if the derivative $\DH$ is not considered as a function of $r$
itself, but only as a {\it single}  derivative. It is in fact immediate to see that in this
case: 
$$ \Phi^\prime [r_\eta]-\big[\DH,\Phi\big] \ = \
 \DH L_{r_\eta}\DH^{-1}+\DH(\DH+L_r)-\DH^2(\DH+L_r)\DH^{-1} \ = \ 0 \ . $$

Unluckily, this fortunate event does not repeat itself in the more difficult case of
hereditariness. Indeed, to prove that operator $\Phi(r)$ is hereditary the complete form
\eqref{mm1} must be used, and a long computation is necessary, together with several 
integration by
parts, to acquire the desired result. In fact, if the difference
$\Phi\Phi^\prime[V]-\Phi^\prime[\Phi V]$ is subdivided into its four terms due to the four terms of
operator $\Phi^\prime $, namely:
$$ \begin{aligned}
 \Phi^\prime _1[V]&= R_V \\
  \Phi^\prime _2[V]&= L_{\DH V}\DH^{-1} \\
  \end{aligned}\hspace{4em}
  \begin{aligned}
   \Phi^\prime _3[V]&= C_rV\DH^{-1} \\
    \Phi^\prime _4[V]&= r_\eta\DH^{-1}C_V\DH^{-1} \\
\end{aligned}   $$
then the four values for the difference $S_j:=\Phi\Phi_j^\prime[V]-\Phi_j^\prime[\Phi
V]$,\quad $j=1,..,4$ are as follows \footnote{In the following terms, the symbol  $L$, which denotes left multiplication, is omitted to simplify the notation.}
$$S_1 \ = \  \left(\begin{array}{c}
 - R_{r V} \\
 - R_{r_{\eta} ( \DH^{-1}V)} \\
 R_{V} \ \DH \\
 r \ R_{V} \\
 r_{\eta} \ \DH^{-1} \ R_{V} \\
\end{array} \right)
\hspace{5em} S_2 \ = \  \left(\begin{array}{c}
 V_\eta \\
 - 2 r_{\eta} \ V \ \DH^{-1} \\
 - r_{\eta\eta} ( \DH^{-1}V) \ \DH^{-1} \\
 r_{\eta} \ \DH^{-1} \ V_\eta \ \DH^{-1} \\
\end{array} \right) $$

$$S_3 \ = \  \left(\begin{array}{c}
 r \ V \\
 - V \ r \\
 r_{\eta} \ V \ \DH^{-1} \\
 - V \ r_{\eta} \ \DH^{-1} \\
 r_{\eta} ( \DH^{-1}V )\ r \ \DH^{-1} \\
 - r \ r_{\eta} \ (\DH^{-1}V) \ \DH^{-1} \\
 r_{\eta} \ \DH^{-1} \ r \ V \ \DH^{-1} \\
 - r_{\eta} \ \DH^{-1} V \ r \ \DH^{-1} \\
\end{array} \right)
\hspace{3em} S_4 \ = \  \left(\begin{array}{c}
 r_{\eta} \ V \ \DH^{-1} \\
 - r_{\eta} \ R_{V} \ \DH^{-1} \\
 r_{\eta} \ \DH^{-1} \ R_{V_\eta} \ \DH^{-1} \\
 r_{\eta} \ \DH^{-1} \ R_{r V} \ \DH^{-1} \\
 r_{\eta \eta} \ \DH^{-1} \ V \ \DH^{-1} \\
 - r_{\eta} \ \DH^{-1} \ V_\eta \ \DH^{-1} \\
 - r_{\eta \eta} \ \DH^{-1} \ R_{V} \ \DH^{-1} \\
 r \ r_{\eta} \ \DH^{-1} \ V \ \DH^{-1} \\
 - r \ r_{\eta} \ \DH^{-1} \ R_{V} \ \DH^{-1} \\
 - r_{\eta} \ \DH^{-1} \ r \ V \ \DH^{-1} \\
 - r_{\eta} \ \DH^{-1} \ r_{\eta} \ (\DH^{-1}V) \ \DH^{-1} \\
 r_{\eta} \ \DH^{-1} \ r_{\eta} \ \DH^{-1} \ V \ \DH^{-1} \\
 - r_{\eta} \ \DH^{-1} \ r_{\eta} \ \DH^{-1} \ R_{V} \ \DH^{-1} \\
 r_{\eta} \ \DH^{-1} \ R_{r_{\eta} ( \DH^{-1}V)} \ \DH^{-1} \\
\end{array} \right) $$

It is clear that even in the automatic procedure the form \eqref{mm1} implies that the
actual result of the term that has to be symmetric in the exchange $W\leftrightarrow V$
is sufficiently long:

$$\left(
  \begin{array}{c}
 r_{\eta} \ \DH^{-1}(\DH^{-1}W) \ V_{\eta} \\
 r_{\eta} \ \DH^{-1}W  \ V  \\
 r_{\eta} \ \DH^{-1}r_{\eta} \ (\DH^{-1}V  \ (\DH^{-1}W)) \\
 r_{\eta} \ \DH^{-1}(\DH^{-1}W) \ r  \ V  \\
 r_{\eta} \ \DH^{-1}(\DH^{-1}W) \ r_{\eta} \ (\DH^{-1}V) \\
 W_{\eta} \ V  \\
 V_{\eta} \ W  \\
 r_{\eta\eta} \ \DH^{-1}V  \ (\DH^{-1}W) \\
  -r_{\eta} \ \DH^{-1}r_{\eta} \ (\DH^{-1}(\DH^{-1}W) \ V ) \\
  -r_{\eta} \ \DH^{-1}V  \  r  \ (\DH^{-1}W) \\
  -r_{\eta} \ \DH^{-1}r_{\eta} \ (\DH^{-1}V) \ (\DH^{-1}W) \\
  -r_{\eta\eta} \ \DH^{-1}(\DH^{-1}W) \ V  \\
  \end{array}
\right) \quad + \quad 
  \left(
  \begin{array}{c}
 r  \ W  \ V  \\
 r  \ V  \ W  \\
 r  \ r_{\eta} \ \DH^{-1}V  \ (\DH^{-1}W) \\
  -r  \ r_{\eta} \ \DH^{-1}(\DH^{-1}W) \ V \\
  -W  \ r  \ V   \\
  -W  \ r_{\eta} \ (\DH^{-1}V) \\
  -V  \ r_{\eta} \ (\DH^{-1}W) \\
  -V  \  r  \  W  \\
  -r_{\eta} \  (\DH^{-1}W) \  V \\
   -r_{\eta\eta}  \ (\DH^{-1}V) \ (\DH^{-1}W) \\
 r_{\eta} \ (\DH^{-1}V) \ r  \ (\DH^{-1}W) \\
  -r  \ r_{\eta} \ (\DH^{-1}V) \ (\DH^{-1}W)
  \end{array}
\right)$$
\medskip

The final step in the automatic computations has been that of proving that this term is
indeed symmetric in the exchange between $W$ and $V$, a fact that as already mentioned
has required several integrations by part, many of which proved to be more
conveniently solved by hand rather than by the automatic procedure.

On the other hand, when the necessary {\it  macros} for the symbolic language are ready
for the non-commutative Burgers' mirror equation, then it is only a matter of care to use
them again with some similar equation. For instance, all the corresponding properties of
the direct non-commutative Burgers' equation: $s_t=\Psi(s)s_x=s_{xx}+2ss_x$ have again
been found with respect to the corresponding {\it new}  derivation: $\DD:=(D+C_s)$.

\section{Appendix }
The aim of this Appendix is twofold; indeed, it collects, in its initial part,  some background notions and
 definitions used throughout the whole article while, in the second part, results on the non-Abelian Burgers,
 in \cite{Carillo:Schiebold:JNMP2012}  are briefly recalled.
\subsection{some background definitions}

\begin{definition}\label{definition1}
(Symmetry)

\noindent
Given an evolution equation $u_t=K(u)$,  where $ u (x,\cdot) \in M $ and  {\hbox{$ K : M \rightarrow TM$}} is
 an appropriate $ C^{\infty} $ vector field on a manifold $ M$, a map $ \sigma : M \rightarrow TM $
is said to be an
\it {infinitesimal symmetry generator }\rm
(for short symmetry) if it leaves the evolution equation itself
 invariant under the infinitesimal transformation $       u \rightarrow u + \varepsilon \sigma$.
\end{definition}
As stated in \cite{Fuchssteiner1979}, if $\sigma$ and $K$ are in involution, i.e. if $[\sigma,K]$ is
identically zero, then $\sigma$ is a symmetry of the given nonlinear evolution equation.

\begin{definition}\label{definition2}
(Strong Symmetry)

\noindent
An operator-valued
function $\Gamma(u)$ is called a strong symmetry of  $u_t=K(u)$ if, for every symmetry $V$
it admits, the vector field $\Gamma(u) V$ is again a symmetry.
\end{definition}

If  $\Gamma$ is a strong symmetry of   $u_t=K(u)$, as proved in \cite{Fuchssteiner1979}, 
the condition
$\Gamma'[K]V=K'[\Gamma V]-\Gamma K'[V]$
holds for any vector field $V$.
\begin{definition}\label{definition3}
(Hereditariness)

\noindent
An an operator-valued function
$\Gamma$ is called hereditary if for every $u\in M$ where $\Gamma$ is defined, the bilinear form
\begin{equation} \label{hereditary}
   \Gamma\ \Gamma'[V]W - \Gamma'[\Gamma V ] W,
\end{equation}
is symmetric in $V$, $W\in T_u  M$.
\end{definition}

\subsection{The non-Abelian Burgers hierarchy}

This Section is devoted to a brief overview on known results concerning the non-Abelian
Burgers equation, the related recursion operator as well as the hierarchy it generates.
Crucial tool is a non-Abelian generalization of the Cole-Hopf transformation connecting
the Burgers equation to the linear heat equation. Given the non-Abelian heat equation
\begin{equation}\label{Heat}
u_t=K(u) ~~, ~~ K(u) = u_{xx}
\end{equation}
and the non-Abelian Burgers equation
\begin{equation}\label{G}
s_t=G(s) ~~, ~~ G(s) = s_{xx}+2s s_x,
\end{equation}
they are connected via the Cole-Hopf transformation  $s=u^{-1}u_x$, which can be written under the form of
 B\"acklund transformation:
\begin{equation}\label{B}
B(u,s)=0~, ~~~ \text{where}~~~ B(u,s)= us-u_x~.
\end{equation}
This connection, given the trivial recursion operator $D$, admitted by the heat
equation, according to  \cite{SIMAI2008, Carillo:Schiebold:JNMP2012}, allows to construct
the recursion operator $\Psi(s)$, admitted by the Burgers equation, that is

\begin{equation} \label{rec RE}
      \Psi(s) = (D+C_s) (D + L_s) (D+C_s)^{-1} ,
\end{equation}
which can also be written as
\begin{equation} \label{rec}
   \Psi(s)=D + L_s + R_{s_x} (D+C_s)^{-1}.
\end{equation}
The latter is the form of the recursion operator also obtained by G\"urses, Karasu and
Turhan \cite{GKT} via a Lax pair    representation of the non-commutative Burgers
hierarchy. Then, the following hierarchies, respectively \eqref{heat hierarchy} and
\eqref{bu hierarchy}, are constructed on application of the trivial recursion operator
$D$, admitted by the heat  equation and the recursion operator $\Psi(s)$ in \eqref{rec
RE}
\begin{equation} \label{heat hierarchy}
   u_{t_n} = D^{n-1}u_x, \quad n\geq 1,
\end{equation}
the lowest members of which read
\begin{eqnarray*}
   u_{t_1} &=& u_x, \\
   u_{t_2} &=& u_{xx} , \\
   u_{t_3} &=& u_{xxx} .
\end{eqnarray*}
and
\begin{equation} \label{bu hierarchy}
   s_{t_n} = \Psi(s)^{n-1}s_x, \quad n\geq 1,
\end{equation}
the lowest members of which read
\begin{eqnarray*}
   s_{t_1} &=& s_x, \\
   s_{t_2} &=& s_{xx} + 2ss_x, \\
   s_{t_3} &=& s_{xxx} + 3ss_{xx} + 3s_x^2 + 3 s^2s_x.
\end{eqnarray*}
The algebraic properties of the operator $\Psi(s)$, firstly obtained  in \cite{SIMAI2008}
and, independently, in \cite{GKT}, are  studied in  \cite{Carillo:Schiebold:JNMP2012}
where $\Psi(s)$ is proved to be a strong symmetry, which is also hereditary.

{\bf Remark} Finally, note that, as expected, if the unknown operator functions $s$ and $r$, 
respectively, in the  non-Abelian Burgers \eqref{bu hierarchy} and mirror non-Abelian 
Burgers hierarchy \eqref{mbu hierarchy}  are replaced by a real valued unknown function 
$v$, then, the commutative Burgers hierarchy is obtained. Furthermore, when $v$ is 
substituted to  $s$ and  $r$, in turn, in the expressions of the two recursion operators 
$\Psi(s)$, in  \eqref{rec}, 
and $\Phi(r)$, in \eqref{mrec2}, they both reduce to the usual (commutative) form of the 
Burgers hereditary recursion operator, that is
\begin{equation} \label{rec-comm}
   \Phi(v)  \equiv \Psi(v) = D + v + {v_x} D^{-1}.
\end{equation}
Hence, the  (commutative) Burgers hierarchy follows as a special case of both the non-Abelian Burgers hierarchies \eqref{bu hierarchy} and  \eqref{mbu hierarchy}.

    \bibliographystyle{model1-num-names}
    \bibliography{sample.bib}

\begin{thebibliography}{00}
\bibitem{Bill2010} D. E. Baldwin, W. A  Hereman, {\sl A symbolic algorithm for computing recursion operators of nonlinear partial differential equations},. Int. J. Comput. Math. {\bf 87}, no. 5, 1094--1119 (2010). 
%

\bibitem{[77]}  S~Carillo,  %
{\sl Nonlinear Evolution Equations: B\"acklund Transformations and B\"acklund Charts},\
Acta Applicandae Mathematicae, {\bf 122}, n.ro 1, pp. 93-106, (2012); 

\bibitem{BS1}
     S.~Carillo and B.~Fuchssteiner,
        \emph{The abundant symmetry structure of hierarchies of nonlinear equations obtained by reciprocal links.}
      J. Math. Phys. {\bf 30}, {1606--1613}  (1989).

 \bibitem{Carillo:Schiebold:JMP2009}
          S.~Carillo and C.~Schiebold,
          \emph{Non-commutative KdV and mKdV hierarchies via recursion methods.}
          J. Math. Phys. {\bf 50}, 073510 (2009).

   \bibitem{Carillo:Schiebold:JMP2011}
          S.~Carillo and C.~Schiebold.
          \emph{Matrix Korteweg-de Vries and modified Korteweg-de Vries hierarchies:
                        Non-commutative soliton solutions.}
          J. Math. Phys. {\bf 52}, 053507 (2011).

\bibitem{SIMAI2008}   S.~Carillo and  C. Schiebold  {\sl A non-commutative operator-hierarchy of Burgers
equations and B\"acklund transformations}, in {\it Applied and Industrial Mathematics in
Italy III: Selected Contributions from the 9th SIMAI Conference}, E. De Bernardis, R.
Spigler e V. Valente Ed.s, SERIES ON ADVANCES IN MATHEMATICS FOR APPLIED SCIENCES,
vol.82, pp. 175 --185, World Scientific Pubbl., Singapore, 2009;

\bibitem{Carillo:Schiebold:JNMP2012}   S.~Carillo and  C. Schiebold {\sl On the recursion operator for
the non-commutative Burgers hierarchy},\ \  { J. Nonlinear Math. Phys.},   {\bf 19} n.ro
1 (2012);

\bibitem{CMS-2015} S.~Carillo, M.~Lo Schiavo and C.~Schiebold, {\sl B\"acklund Transformations \!and\! Non Abelian Nonlinear Evolution Equations: a novel B\"acklund Chart}, submitted, (2015);

   \bibitem{Carl Schiebold} B.~Carl and C.~Schiebold.
      \emph{Nonlinear equations in soliton physics and operator ideals}.
      Nonlinearity {\bf 12}, 333--364 (1999).

   \bibitem{Carl Schiebold 1} B.~Carl and C.~Schiebold.
      \emph{Ein direkter Ansatz zur Untersuchung von Solitonengleichungen}.
      Jahresber. Deutsch. Math.-Verein. {\bf 102}, 102--148 (2000).

\bibitem{Cole:1951}
{J.D.~Cole}
 \emph{On a quasilinear parabolic equation occuring in aerodynamics}
{ Quart.App. Math.}{\bf 92}, {25--236} {(1951).}

\bibitem{FokasFuchssteiner:1981}
      A.S.~Fokas and B.~Fuchssteiner.
      \emph{B\"acklund transformation for hereditary symmetries.}
      Nonlin. Anal., Theory Methods Appl. {\bf 5}, No. 4, 423--432 (1981).

    \bibitem{Fuchssteiner:Carillo:1989a}
    {B.~Fuchssteiner and S.~Carillo}
    {\it Soliton structure versus singularity analysis: Third order completely integrable nonlinear equations in 1+1 dimensions.}
    {Physica} {A 152}, pp. {467-510}, {(1989)}.

   \bibitem{Fuchssteiner1979} 
          B.~Fuchssteiner.
          {Application of hereditary symmetries to nonlinear evolution equations.}
          {\textit{Nonlin. Anal., Theory Methods Appl. }}{\bf 3}, No. 6, 849--862 (1979).

\bibitem{FoFu2}B.~Fuchssteiner and A.S.~Fokas:
Symplectic Structures, Their B\"acklund Transformations and
Hereditary Symmetries,
Physica vol. 4D, pp. 47-66, 1981.

  \bibitem{FuOWi}
      \newblock B. Fuchssteiner, W. Oevel, and W. Wiwianka.
      \newblock \emph{Computer-algebra methods for investigating hereditary operators
                      of higher order soliton equations.}
      \newblock Computer Phys. Commun. {\bf 44}, 47--55 (1987).

   \bibitem{GKS}
      \newblock M.~G\"urses, A.~Karasu, and V.V.~Sokolov.
      \newblock \emph{On construction of recursion operators from Lax representation.}
      \newblock J. Math. Phys. {\bf 40}, 6473--6490 (1999).

   \bibitem{GKT}
      \newblock M.~G\"urses, A.~Karasu, and R.~Turhan.
      \newblock \emph{On non-commutative integrable Burgers equations.}
      \newblock J. Nonlinear Math. Phys. {\bf 17}, 1--6 (2010).

\bibitem{Hopf:1950}
{E.~Hopf}
 \emph{The partial differential equation $ u_{t} + u u_{x} = m
uu_{xx}$.}
{Comm. Pure Appl. Math.}
{\bf 3}, {201--230} (1950).

   \bibitem{Ku}
      \newblock B.A.~Kupershmidt.
      \newblock \emph{On a group of automorphisms of the noncommutative Burgers hierarchy.}
      \newblock J. Nonlinear Math. Phys. {\bf 12}, No. 4, 539--549 (2005).

   \bibitem{LRB}
          D.~Levi, O.~Ragnisco  and M.~Bruschi.
          \emph{Continuous and discrete matrix Burgers' hierarchies.}
          Il Nuovo Cimento {\bf 74B}, 33--51 (1983).

\bibitem{Olver}
    P.J.~Olver.
    {\it Evolution equations possessing infinitely many symmetries.}
    J.~Math.~Phys.~{\bf 18}, 1212-1215 (1977).

   \bibitem{Olver:Sokolov}
          P.J.~Olver and V.V.~Sokolov.
          \emph{Integrable evolution equations on nonassociative algebras.}
          Comm. Math. Phys. {\bf 193}, 245--268 (1998).

    \bibitem{Rogers:Carillo:1987b}C.~Rogers and S.~Carillo \ {\sl On
    reciprocal properties  of  the
    Caudrey-Dodd-Gibbon  and
     Kaup-Kupershmidt hierarchies}, \
     Physica Scripta, {\bf 36}, (1987), 865-869.

    \bibitem{Schiebold 1} C.~Schiebold.
      \emph{Explicit solution formulas for the matrix-KP}.
      Glasgow Math. J. {\bf 51}, 147--155 (2009).

        \bibitem{Schiebold-6dic1} C.~Schiebold.
     {\it Cauchy-type determinants and integrable systems}.
Linear Algebra and its Applications {\bf 433}, 447�-475 (2010).

  \bibitem{Schiebold-6dic2} C.~Schiebold.
     {\it Noncommutative AKNS systems and multisoliton solutions to the matrix
sine-Gordon equation},
Discr. Cont. Dyn. Systems Suppl. {\bf 2009}, 678�-690 (2009).

  \bibitem{Schiebold-6dic3} C.~Schiebold.
     {\it The noncommutative AKNS system: projection to matrix systems, countable
superposition and soliton-like solutions},
J. Phys. A {\bf 43}, 434030 (2010).

   \bibitem{Schiebold2010}
      \newblock C.~Schiebold.
      \newblock \emph{Structural properties of the noncommutative KdV recursion operator.}
      \newblock J. Math. Phys. {\bf 52}, 113504 (2011).

\end{thebibliography}

\subsection*{Acknowledgements}

The financial support of G.N.F.M.-I.N.d.A.M.,  I.N.F.N. and \textsc{Sapienza}  University of Rome, Italy are gratefully acknowledged.
C. Schiebold wishes also to thank S.B.A.I. Dept.  and \textsc{Sapienza} University of Rome for the kind hospitality.

\end{document}